\documentclass[11pt]{article}

\usepackage{amssymb,amscd,array}

\let\ssection=\section
\renewcommand{\section}{\setcounter{equation}{0}\ssection}

\newcommand{\bbR}{\mathbb{R}}

\newcommand{\bbC}{\mathbb{C}}

\newcommand{\Diff}{\mathrm{Diff}}

\newcommand{\cF}{{\mathcal{F}}}
\newcommand{\cD}{{\mathcal{D}}}

\newcommand{\Hom}{\mathrm{Hom}}

\newcommand{\PSL}{\mathrm{PSL}}
\newcommand{\SL}{\mathrm{SL}}
\newcommand{\Sl}{\mathrm{sl}}

\newcommand{\Vect}{\mathrm{Vect}}

\newcommand{\Hol}{\mathrm{Hol}}

\begin{document}



\def\d{\delta}
\def\S{\Sigma}
\def\n{\nabla}
\def\G{\Gamma}
\def\g{\gamma}
\def\om{\omega}
\def\r{\rho}
\def\a{\alpha}
\def\b{\beta}
\def\s{\sigma}
\def\vfi{\varphi}
\def\l{\lambda}
\def\m{\mu}
\def\implies{\Rightarrow}

\oddsidemargin .1truein
\newtheorem{thm}{Theorem}[section]
\newtheorem{lem}[thm]{Lemma}
\newtheorem{cor}[thm]{Corollary}
\newtheorem{pro}[thm]{Proposition}
\newtheorem{ex}[thm]{Example}
\newtheorem{rmk}[thm]{Remark}
\newtheorem{defi}[thm]{Definition}

\title{Projectively Invariant Cocycles of Holomorphic 
Vector Fields on an Open Riemann Surface}
\date{}
\maketitle
\noindent S. BOUARROUDJ \\
{\small Department of Mathematics, Keio University, Faculty of Science \& Technology. 
3-14-1, Hiyoshi, Kohoku-ku, Yokohama, 223-8522, Japan.  E-mail:sofbou@math.keio.ac.jp}\\

\noindent H. GARGOUBI \\
{\small I.P.E.M., route de Kairouan, 5019 Monastir, Tunisia. E-mail:hichem.gargoubi@ipeim.rnu.tn}

\begin{abstract} Let $\Sigma$ be an open Riemann surface and $\Hol (\Sigma )$ be the 
Lie algebra of holomorphic vector fields on $\Sigma.$ We fix a projective structure 
(i.e. a local $\SL_2(\bbC)-$structure) on $\Sigma.$ 
We calculate the first group of cohomology  
of $\Hol(\Sigma )$ with coefficients in the space of linear holomorphic operators 
acting on tensor densities, vanishing on the Lie algebra $\Sl_2(\bbC)$. 
The result is independent on the  
choice of the projective structure. We give explicit formul{\ae} of 1-cocycles 
generating this cohomology group.
\end{abstract}
\section{Introduction}
The first group of cohomology of the Lie algebra of (formal) vector fields on the circle 
$S^1$ with  coefficients in the space $\Hom ({\cal F}_\l,{\cal F}_{\mu}),$ where 
${\cal F}_{\l}$ is the space of tensor densities of degree $\l$ on $S^1,$ was first 
calculated in \cite{gf}. This group of cohomology measures all extensions of exact 
sequences  $0\rightarrow {\cal F}_\mu \rightarrow \cdot 
\rightarrow {\cal F}_\l\rightarrow 0$ of modules. The first group of cohomology of the 
Lie algebra $\Vect (S^1)$ of (smooth) vectors fields on the circle  with  
coefficients in the space of linear differential  
operators acting from $\l-$densities to $\mu-$densities, 
vanishing on the subalgebra
$\Sl_2(\bbR)\subset \Vect(S^1)$, was calculated in \cite{bo}.  
This group of cohomology appears as an obstruction to the equivariant quantization 
(see~\cite {bo}, \cite{ga}). The computation is based on the following observation: 
any 1-cocycle vanishing on the subalgebra $\Sl_2(\bbR)$ is an $\Sl_2(\bbR)-$invariant 
operator. The $\Sl_2(\bbR)-$invariant differential operators acting on tensor densities, which
are called ``Transvectants'', were classified by Gordan (see \cite{g,jp}). To find the 
1-cocycles generating the group of cohomology means, therefore, to determine which 
from the Transvectants are  1-cocycles. 

In this paper, we study the complex analog of the above group of 
cohomology on an open Riemann surface $\Sigma$ endowed with a flat projective 
structure.

The aim of this paper is to describe the first group of cohomology 
\begin{equation}
\label{gr}
H^1(\Hol (\Sigma), \Sl_2 (\bbC); \cD_{\l,\mu}),
\end{equation} 
of holomorphic vector fields on $\Sigma$ with coefficients in the space of linear 
differential operators acting on tensor densities, vanishing on the Lie algebra 
$\Sl_2(\bbC)$. We give explicit  
formul{\ae} of 1-cocycles generating the group (\ref{gr}). These 1-cocycles are 
the complex analog of the 1-cocycles given in \cite{bo}. 

The main tool of this paper is the existence of affine and projective connection on 
any open Riemann surface (see \cite{gu,kh}). These notions has been recently used  
in \cite{wa} to compute the second group of cohomology of the Lie algebra 
$\Hol (\Sigma)$ with coefficients in the space of $\l-$densities.
\section{Affine and Projective Structure}
Let $\S$ be a Riemann surface, and let $\{ U_\a, z_\a\}$ be an atlas of $\Sigma$. 

A holomorphic affine connection is a family of holomorphic functions $\G(z_\a)$ on 
$U_\a$ such that for non-empty $U_\a\cap U_\b,$ we have 
$$
\G (z_\b)\,\, \frac{d z_\b}{d z_\a}=\G (z_\a) + \frac{d^2 z_\b}
{d z_\a^2}\frac{d z_\a}{d z_\b}\cdot
$$
Affine connection exists in any open Riemann surface in contrast with the compact 
case where affine connection exists only if the genus of $\Sigma$ is one    
(see \cite{gu}).

A holomorphic projective connection is a family of holomorphic functions $R(z_\a)$ on 
$U_\a$ such that for non-empty $U_\a\cap U_\b,$ we have 
$$
R(z_\b) \,\,\left (\frac{dz_b}{dz_\a}\right )^2=R(z_\a) + S(z_\b,z_\a),
$$
where $S(z_\b,z_\a)=\left (  \frac{d^2 z_\b}{d z_\a^2}\right )'-\frac{1}{2}
\left (  \frac{d^2 z_\b}{d z_\a^2}\right )^2$ is the Schwarzian derivative.\\

Recall that any holomorphic affine connection $\G$ defines naturally a holomorphic 
projective connection given in a local coordinates $z$ by  
\begin{equation}
\label{connetion}
R(z)=\frac{d \G(z)}{dz}-~\frac{1}{2}\G^2(z).
\end{equation}
For any projective connection $R$ there exits locally an affine connection 
$\Gamma$ satisfying (\ref{connetion}).

We say that two affine connections are projectively equivalent if they define the  
same projective connection.

Let us define the notion of projective structure.

A Riemann surface admits a projective structure if there exists an atlas of charts 
$\{U_\a, z_\a\}$ such that the coordinate change $z_\a\circ z_\b^{-1}$ are projective 
transformations. 

In this case, the Lie algebra $\Sl_2 (\bbC),$ in each chart of the projective structure,  
is generated by the following vector fields
\begin{equation}
\label{champ}
\frac{d}{dz},\, z\frac{d}{dz},\, z^2\frac{d}{dz}.
\end{equation} 
There exists a 1-1 correspondence  between  projective structures and projective connections 
on an open Riemann surface (cf. \cite{gu}). We use implicitly this correspondence along this 
paper.   
\section{$\Hol(\Sigma)$-module Structures on the Space of Differential\break 
Operators}
The modules of linear differential operators on the space of tensor densities on a 
(real) smooth manifold has been studied in series of recent papers
(see \cite{b, bo, do, do2, ga, go, lmt, lo}). Note that this space viewed as 
a module over the Lie algebra of vector fields has already been studied in 
the classical monograph \cite{w}.  

Let us give the definition of the natural two-parameter family of modules
over the Lie algebra of holomorphic vector fields on the space of linear 
differential operators.

\vskip 0,3cm
\noindent
\subsection{ Tensor Densities}
Let $\S$ be an open Riemann surface. Fix an affine connection $\Gamma$ on it.

Space of tensor densities on $\S$, noted ${\cal F}_{\l}$, is the space of 
sections of the line bundle $ (T^*\S)^{\otimes \l}$, where 
$\l\in \bbC.$ This bundle is of course trivial, since any (holomorphic) bundle on an open      
Riemann surface is holomorphically trivial. 

Fix a global section $dz^\l$ on $\cF_\l.$ Any $\l-$density can be written 
in the form $\phi\,\, dz^\l.$ Let us recall the definition of a covariant 
derivative of tensor densities. Let $\n$ be the covariant derivative associated 
to  the affine connection $\G.$   If $\phi \in \cal F_\l,$ then 
$\n\,\phi\in \Omega^1(\S)\otimes {\cal F}_\l$ given by the formula
$$ \n\, \phi=\frac{d \phi}{dz}-\l\G \phi.$$

The standard action of $\Hol(\Sigma)$ on $\cal F_\l$ reads as follows (cf. \cite{wa}):
\begin{equation}
\label{dens}
\nonumber
L_{X}^{\l}(\phi)=X\n \phi +\l \phi \,\n X, 
\end{equation}
where $X=X(z) \displaystyle \frac{d}{dz} \in \Hol (\S).$
\subsection{$\Hol (\Sigma)-$Module of Differential Operators} 
Consider differential operators acting on tensor densities:
\begin{equation}
\label{Op}
A:{\cal F}_{\lambda}\to{\cal F}_{\mu}.
\end{equation}
In local coordinates $z$, any operator $A$ can be written in the form
$$
A=a_k(z)\frac{d^k}{dz^k}+\cdots +a_0(z),
$$
where $a_i$, for $i=1,k$, are holomorphic functions on $z.$\\
\noindent A two-parameter family of actions of $\Hol(\Sigma) $ on the 
space of differential operators is defined by
\begin{equation}
\label{act}
L_X^{\l,\mu}(A)=L_X^\mu \circ A-A\circ L_X^\l.
\end{equation}
Denote by ${\cal D}_{\lambda,\mu}(\Sigma)$ the space of operators (\ref{Op}) endowed 
with the defined $\Hol(\Sigma)$-module structure (\ref{act}).

\section{Main Result}
Assume $\Sigma$ endowed with a projective structure, this defines locally 
an action of the Lie group $\SL_2(\bbC)$ on $\Sigma.$ 

Consider cochains on $\Hol(\Sigma)$ with values in $\cD_{\l,\mu},$ vanishing on the 
Lie algebra $\Sl_2(\bbC).$ One, therefore, obtains the so-called relative cohomology 
of the Lie algebra $\Hol(\Sigma)$, namely 
$$
H^1(\Hol(\Sigma),\Sl_2(\bbC);\cD_{\l,\mu} (\Sigma )),
$$ 
(see \cite{f}).

The purpose of this paper is the following:
\begin{thm}
\label{pre}
The first group of cohomology $H^1(\Hol(\Sigma), \Sl_2(\bbC);{\cal D}_{\l,\mu}(\Sigma ))$ is 
one-dimensional in the following cases:\\
(a) $\mu-\l=2, \l\not=-1/2,$\\
(b) $\mu-\l=3, \l\not=-1,$\\
(c) $\mu-\l=4, \l\not=-1/2,$\\
(d) $(\l,\mu)=(-4,1),(0,5)$.\\
Otherwise, this cohomology group is trivial.
\end{thm}
This Theorem generalizes the result of \cite{bo} in the case of the circle $S^1.$
 
Note that this result does not depend on the choice of the projective structure.
\section{Construction of the 1-cocycles}
In this section, we give explicit formula for the 1-cocycles generating the 
nontrivial cohomology classes from Theorem \ref{pre}. Given a projective structure 
on $\Sigma,$ we prove that there is a canonical choice of the 1-cocycles vanishing 
on $\Sl_2(\bbC).$

Fix (locally) an affine connection $\Gamma$ related to the projective structure. 
Denote by $R$ the projective connection associated to $\G$ (see section 2.)
\begin{lem}
\label{ope}
The following linear differential operators 
$$
\begin{array}{ccl}
{\cal I}_3&=&R,\\
{\cal I}_4&=&R\,\n\, -\displaystyle \frac{\l}{2} \n R,\\[3mm]
{\cal I}_5&=&R\,\n^2\,-\displaystyle \frac{2\l+1}{2}\n R\, \n\,+\frac{\l(2\l+1)}{10}\n^2R+
\frac{\l(\l+3)}{5}R^2,\\[3mm]
{\cal I}_6&=&R\,\n^3\,-\displaystyle \frac{3}{2}\n R\, \n^2\,+
\left (\frac{3}{10}\n^2R+\frac{4}{5}R^2\right )\n\,,\\[3mm]
{\cal I}_6'&=&R\,\n^3\,+\displaystyle \frac{9}{2}\n R\, \n^2\,+
\left (\frac{63}{10}\n^2R+\frac{4}{5}R^2\right )\n\,+
\frac{14}{5}\n^3R+\frac{8}{5}R\n\, R,

\end{array}
$$
are globally defined in ${\cal D}_{\l,\l+2}(\Sigma ),$ 
${\cal D}_{\l,\l+3}(\Sigma ),$ ${\cal D}_{\l,\l+4}(\Sigma ),$ ${\cal D}_{0,5}(\Sigma ),$ 
${\cal D}_{-4,1}(\Sigma ),$ respectively, and depend only on the projective class of the 
connection $\Gamma.$

\end{lem}
{\bf Proof.} Since the surface $\Sigma$ is projectively flat, the connection  $R$ defines a 
2-density on $\Sigma.$ Then the operators of the above Lemma are globally defined. 
Let us prove that the operators depend only on the projective class of the 
connection $\Gamma.$

Let ${\cal I}_4=R\,\n\,+\alpha  \nabla R .$ Denote by 
$\tilde {\cal I}_4$ the operator ${\cal I}_4$ written with respect to a connection 
$\tilde \Gamma$ which is projectively equivalent to $\Gamma.$ After an easy calculation one 
has 
$$
{\cal I}_4=\tilde {\cal I}_4+(2\alpha+\l)(\tilde \Gamma-\Gamma)R.
$$  
Then ${\cal I}_4=\tilde {\cal I}_4$ if and only if $\alpha=-\l/2.$ The proof is 
analogous for 
the operators $ {\cal I}_5, {\cal I}_6$  and ${\cal I}_6'.$ 
 
\begin{thm}
\label{deu}
(i) For every $\l,$ there exist unique (up to constant) 1-cocycles 
$$
\begin{array}{c}
{\cal J}_3:\Hol (\S)\rightarrow {\cal D}_{\l,\l+2}\\
{\cal J}_4:\Hol (\S)\rightarrow {\cal D}_{\l,\l+3}\\
{\cal J}_5:\Hol (\S)\rightarrow {\cal D}_{\l,\l+4}\\
\end{array}
$$
vanishing on $\Sl_2(\bbC)$. They are given by the formul\ae :
$$
\matrix{
\displaystyle {\cal J}_3(X)=
\n^3 X -L^{\l,\l+2}_X( {\cal I}_3),\hfill\cr
\noalign{\smallskip}
\displaystyle {\cal J}_4(X)=
\n^3 X\n\, -\frac{\lambda}{2}\n^4 X-L^{\l,\l+3}_X({\cal I}_4), \hfill\cr
\noalign{\smallskip}
\displaystyle {\cal J}_5(X)=
\n^3 X\n^2\,-\frac{2\lambda+1}{2}\n^4X \n \,+
\frac{\lambda(2\lambda+1)}{10}\n^5X -L^{\l,\l+4}_X({\cal I}_5). \hfill\cr
}
$$
For $(\l,\mu)=(0,5), (-4,1)$, respectively, the 1-cocycles vanishing on $\Sl_2(\bbC)$ 
are given by
$$
\matrix{
\displaystyle {\cal J}_6^{0}(X)=&
\n^3 X\,\n^3\, -\displaystyle \frac{3}{2}\n^4 X\, \n^2 \,+
\frac{3}{10}\n^5 X\n \,  -L_X^{0,5}({\cal I}_6),\hfill \cr
\noalign{\smallskip}
\displaystyle {\cal J}_6^{-4}(X)=&\n^3 X\,\n^3\, +\displaystyle 
\frac{9}{2}\n^4 X\, \n^2\,+\frac{63}{10}\n^5 X\n\,+
\frac{14}{5}\n^6X -L_X^{-4,1}({\cal I}_6').
\hfill\cr
}
$$
(ii) The 1-cocycles ${\cal J}_3$,  ${\cal J}_4$ and ${\cal J}_5$ are nontrivial for every $\l$ 
except $\l=-1/2, \l=-1$ and $\l=-3/2$, respectively. The 1-cocycles  ${\cal J}_6^0$ and 
${\cal J}_6^{-4}$ are nontrivial.

\noindent (iii) These 1-cocycles are independent on the choice of the projective structure.
\end{thm} 
\section{$\Sl_2({\small \bbR})-$invariant operators on $S^1$}
\label{rappel}
For almost all $\l$ and $\mu,$ there exits unique (up to constant) $\Sl_2(\bbR)-$invariant 
bilinear differential operators $J_m^{\l,\mu}: \cF_\l\otimes 
\cF_\mu \rightarrow \cF_{\l+\mu+m}$ given by   
\begin{equation}
\label{trans}
J_m^{\lambda,\mu}(\phi, \psi )
=
\sum_{i+j=m} (-1)^i m!
{2\lambda+m-1 \choose i} {2\mu+m-1 \choose j}
{\phi}^{(i)} {\psi}^{(j)},
\end{equation}
called ``Transvectants'' (see \cite{go, jp}). 

Let us recall the results of \cite{bo,go}. \\
The first group of cohomology $H^1(\Vect(S^1), \Sl_2(\bbR);{\cal D}_{\l,\mu}(S^1))$ is 
one-dimensional in the following cases:\\
(a) $\mu-\l=2, \l\not=-1/2,$\\
(b) $\mu-\l=3, \l\not=-1,$\\
(c) $\mu-\l=4, \l\not=-1/2,$\\
(d) $(\l,\mu)=(-4,1),(0,5).$\\
Otherwise, this cohomology group is trivial (see \cite{bo}).
\newpage
This group of cohomology is generated by the following 1-cocycles which are particular 
cases of the Transvectants (\ref{trans}):
\begin{equation}
\label{tra}
\matrix{
\displaystyle J_3(X,\phi)=
X'''\phi,\hfill\cr
\noalign{\smallskip}
\displaystyle J_4(X,\phi)=
X'''\phi'-\frac{\lambda}{2}X^{IV}\phi,\hfill\cr
\noalign{\smallskip}
\displaystyle J_5(X,\phi)=
X'''\phi''-\frac{2\lambda+1}{2}X^{IV}\phi'+
\frac{\lambda(2\lambda+1)}{10}X^{V}\phi,\hfill\cr
\noalign{\smallskip}
\displaystyle J_6(X,\phi)=
X'''\phi'''
-\frac{3}{2}X^{IV}\phi''
+\frac{3}{10}X^{V}\phi' ,\hfill\cr
\noalign{\smallskip}
\displaystyle J_6'(X,\phi)=
X'''\phi''' +\frac{9}{2}X^{IV}\phi''+\frac{63}{10}X^{V}\phi'+
\frac{14}{5}X^{VI}\phi.\hfill\cr
}
\end{equation}
(see \cite{bo, go}).

The 1-cocycles given in Theorem \ref{deu} are the complex analogue of the 1-cocycles (\ref{tra}).  
\section{Proof of the main Theorems}
In this section, we will prove Theorem \ref{pre} and Theorem \ref{deu}.
\subsection{Proof of the Theorem \ref{deu}}
To prove that the operators ${\cal J}_k$, for $k=3,4,5,6,$ are 1-cocycles one has to check  
the 1-cocycle relation. It reads as follows
\begin{equation}
\label{cocy}
{\cal J}_k[X,Y]-L_X^{\l, \l +k-1}({\cal J}_k(Y))+L_Y^{\l, \l +k-1}({\cal J}_k(X))=0,
\end{equation}
where $X, Y\in \Hol(\Sigma).$\\ 
Let us verify (\ref{cocy}) for the operator ${\cal J}_3.$ It is obvious that 
$L_X^{\l,\l+2}({\cal R})$ is a 1-cocycle. It suffices then to verify the relation 
(\ref{cocy}) for the $0-$order  operator $\n^3 X:$
$$
\small{
\matrix{ 
\n^3([X,Y])\phi-L_X^{\l,\l+2}(\n^3 Y \phi)+L_Y^{\l,\l+2}(\n^3 X \phi)=&
(2\n X\n^3 Y+X \n^4 Y -2\n Y \n^3 X-Y \n^4 X)\,\phi \hfill\cr
\noalign{\smallskip}
&-X\n (\n^3 \phi )-(\l+2)\n X \n^3 Y \phi \hfill\cr
\noalign{\smallskip}
&+\n^3 Y(X\n \phi +\l \n X \phi)+Y\n(\n^3 X\phi )\hfill\cr
\noalign{\smallskip}
&+(\l+2)\n Y \n^3 X \phi -\n^3 X(Y\n \phi +\l \n Y \phi)\hfill\cr
\noalign{\smallskip}
\hfill  =&0.\hfill\cr
}
}
$$
Let us prove that the 1-cocycle ${\cal J}_3$ vanishes on $\Sl_2(\bbC).$ Let $X$ be one 
of the vector fields (\ref{champ}). After calculation one has $\n^3 X=2R\n X +X \n R.$ 
It is easy to see that $L_X^{\l, \l+2}({\cal I}_3)=2R\n X +X \n R.$ Hence, one obtains 
${\cal J}_3(X)=0.$  

In the same manner we prove that the operators ${\cal J}_k,$ for $k=4,5,6,$ are 
1-cocycles vanishing on $\Sl_2(\bbC).$ Theorem \ref{deu} (i) is proven.

Let us prove the non-triviality of the 1-cocycle ${\cal J}_3$ for $\l \not =-1/2.$ 
Suppose that the 1-cocycle ${\cal J}_3(\Sigma )$ is trivial, then there exists an operator 
$A\in \cD_{\l,\l+2}$ such that
\begin{equation}
\label{expr}
{\cal J}_3(X)=L_X^{\l,\l+2}(A).
\end{equation}
In a neighborhood of a point $z$, one can choice a local coordinates such that the connection  
$R=0.$ In these coordinates, the 1-cocycle ${\cal J}_3$ coincides with the 1-cocyle $J_3$ of 
(\ref{tra}). Hence the relation (\ref{expr}) implies that the 1-cocycle $J_3$ is trivial which 
is absurd (see section \ref{rappel}).\\
For $\l =-1/2,$ take $A=-2 \n^2 \in \, \cD_{-\frac{1}{2},\frac{3}{2}}(\Sigma ).$ 
One can easily check that  ${\cal J}_3(X)=L_X^{-1/2,3/2}(A).$\\
With the same arguments we prove the non-triviality of the 1-cocyles ${\cal J}_k,$ for $k=4,5,6.$ 
Theorem \ref{deu} (ii) is proven.  
\subsection{Proof of Theorem \ref{pre}}
Let us prove that the dimension of the group of cohomology $H^1(\Hol (\S), \Sl_2(\bbC); 
\cD_{\l,\mu}(\Sigma ))$ is bounded by the dimension of the group of cohomology 
$H^1(\Vect (S^1), \Sl_2(\bbR); \cD_{\l,\mu}(S^1)).$ Let 
$C$ and $C'$ be two 1-cocycles in $H^1(\Hol (\S), \Sl_2(\bbC); \cD_{\l,\mu}(\Sigma )).$ 
We will prove that $C$ and $C'$ are cohomologous. Denote  $\tilde C$ and $\tilde C'$ the 
restriction of $C$ and $C'$ on a neighborhood of a point of $\S.$ The operators 
$\tilde C$ and $\tilde C'$ define 1-cocycle in 
$H^1(\Vect (S^1), \Sl_2(\bbR); \cD_{\l,\mu}(S^1)).$ These 1-cocycles are equal (up to 
constant); since the unique $\Sl_2(\bbR)-$invariant linear differential operators are 
given as in (\ref{tra}). It follows that the 1-cocyle $C$ and $C'$ are cohomologous. Now from 
the construction of the 1-cocyles given in Theorem (\ref{deu}) follows Theorem (\ref{pre}). 
\section{Final Remark}
The group of cohomology $H^1(\Diff (S^1),\PSL_2(\bbR);\cD_{\l,\mu})$ 
is one-dimension, for generic $\l,$ generated by the following 1-cocyles
\begin{equation}
\label{sch}
\small{
\matrix{
{\cal S}_{\lambda}(f)=&S(f),
\hfill\cr\noalign{\bigskip}
\displaystyle{\cal T}_{\lambda}(f)=&S(f)\displaystyle \frac{d}{dx}-
\frac{\lambda}{2} S(f)^{\prime}, 
\hfill\cr\noalign{\bigskip}
{\cal U}_{\lambda}(f)=&\displaystyle S(f)\frac{d^2}{dx^2}-
\frac{2\lambda+1}{2} S(f)^{\prime}\frac{d}{dx}
+\frac{\lambda(2\lambda+1)}{10} S(f)^{\prime\prime}-
\frac{\lambda(\lambda+3)}{5} S(f)^2,\hfill\cr\noalign{\bigskip}
{\cal V}_0(f)=&\displaystyle S(f)\frac{d^3}{dx^3}-
\displaystyle \frac{3}{2} S(f)^{\prime}\frac{d^2}{dx^2}+
\left (\frac{3}{10}S (f)^{\prime\prime}+\frac{4}{5}S(f)^2\right )\frac{d}{dx}
\hfill\cr\noalign{\bigskip}
{\cal V}_{-4} =&S(f)\displaystyle \frac{d^3}{dx^3}+\displaystyle \frac{9}{2} S(f)^{\prime}
\, \frac{d^2}{dx^2}+
\left (\frac{63}{10}S(f)^{\prime\prime}+\frac{4}{5}S(f)^2\right )\frac{d}{dx}+
\frac{14}{5}S(f)^{\prime\prime\prime}+\frac{8}{5}S(f)S(f)^{\prime},
}
}
\end{equation}
where $S(f)$ is the Schwarzian derivative (see \cite{bo}).
It is a remarkable fact to see that the 1-cocycles  (\ref{sch}) on the group $\Diff(S^1)$ have   
the same expression (up to change of sign) than the operators of the Lemma \ref{ope} if one 
replaces the connection $R$ by $-S(f).$ 

\noindent Since the group of biholomorphic maps on a Riemann surface 
is finite dimension (see \cite{kf}), this group does not integrate the Lie algebra of holomorphic vector 
fields. In some sense, the cohomology group (\ref{gr}) contains informations coming from the 
cohomology of the diffeomorphisms of $S^1$ and the cohomology of the Lie algebra of vector 
fields on $S^1.$

\bigskip

{\it Acknowledgments}. It is a pleasure to acknowledge numerous fruitful discussions with 
Ch. Duval, V. Ovsienko and F. Wagemann. The first author is  grateful to the JSPS for the research 
support, and Prof. Y. Maeda and Keio University for their hospitality. 

\vskip 1cm




\end{document}